\documentclass{article}

\usepackage{amssymb} 
\usepackage{amssymb} 
\usepackage{latexsym}
\usepackage{amsmath} 

\newcommand{\W}{{\mathcal W}}
\def\a{a}
\def\bb{b}
\def\n{\nabla}
\def\N{\nabla}

\def\nn{\noindent}

\def\lfleche{\smash{\mathop{\vbox{\hbox to 8mm {{\rightarrowfill}}}}
    \limits^{}_{}}}

\newtheorem{theorem}{Theorem}[section]

\newtheorem{definition}[theorem]{Definition}

\title{On the problem of linearizability of a 3-web}

\author{Zolt\'an Muzsnay}

\begin{document}

\maketitle

\begin{abstract}  
  In this paper we study the linearizability problem for $3$-webs on a
  $2$-dimensional manifold.  With an explicit computation based on the
  theory developed in \cite{GMS}, we examine a 3-web whose
  linearizability was claimed in \cite{GMS}. We show that, contrary to
  the statement of \cite{GL} and \cite{GL_cras}, this particular web
  is linearizable.  We compute explicitly the affine deformation
  tensor and the corresponding flat linear connection adapted to the
  web which linearizes the 3-web.
\end{abstract} 

\bigskip

\begin{description}
\item[AMS Classification:] 53A60, 53C36
  
\item[Keywords:] Webs, linearization, affine structures, Gronwall
  conjecture.
\end{description}

\section{Introduction}

On a two-dimensional real or complex differentiable manifold $M$ a
$3$-web is given by three foliations of smooth curves in general
position.  Two webs ${\mathcal W}$ and $\tilde{{\mathcal W}}$ are
locally equivalent at $p\in M$, if there exists a local diffeomorphism
on a neighborhood of $p$ which exchanges them.  A $3$-web is called
{\it linear} (resp.~{\it parallel}) if it is given by 3 foliations of
straight lines (resp. of parallel lines). A $3$-web which equivalent
to a linear (resp. parallel) web is called linearizable (resp.
parallelizable). An elegant characterization of parallelizable webs
can be given in terms of the Chern connection associated to a $3$-web:
a $3$-web is parallelizable if and only if the curvature of the Chern
connection vanishes.

The problem of finding linearizability criterion is a very natural
one.  Bol suggested a method in \cite {Bol1} how to find a criterion
of linearizability, although he was unable to carry out the
computation.  He showed that the number of projectively different
linear $3$-webs in the plane which are equivalent to a
non-parallelizable $3$-web is finite and less that 17.  The
formulation of the linearizability problem in terms of the Chern
connection was suggested by Akivis in a lecture given in Moscow in
1973. In his approach the linearizability problem is reduced to the
solvability of a system of partial differential equations on the
components of the affine deformation tensor.  Using Akivis' idea
Goldberg determined in \cite{Gol2} the first integrability conditions
of the system.

By using this approach Grifone, Muzsnay and Saab solved the
linearizability problem \cite{GMS}.  They showed that, in the
non-parallelizable case, there exists an algebraic submanifold
${\mathcal A}$ of the space of vector valued symmetric tensors
($S^2T^*\otimes T$) on a neighborhood of $p$, expressed in terms of
the curvature of the Chern connection and its covariant derivatives up
to order $6$, so that the affine deformation tensor is a section of
$S^2T^*\otimes T$ with values in ${\mathcal A}$.  In particular:
\begin{enumerate}
\item The web is linearizable if and only if ${\mathcal A}\neq
  \emptyset$;
\item There exists at most $15$ projectively nonequivalent
  linearizations of a nonparallelizable 3-web.
\end{enumerate}
The expressions of the polynomials and their coefficients which define
$\mathcal{A}$ can be found in \cite{GMS_arxiv}. The criteria of 
linearizability gives the possibility to make explicit computation on
concrete examples to decide whether or not they are linearizable.

Recently Goldberg and Lychagin find similar results on the
linearizability in \cite{GL}, but their method is different from that
of \cite{GMS}.  Despite of the fact, that the two theories concern the
same problem and the final results are very similar, by testing them
on an explicit example they lead to different answers. Indeed,
considering the 3-web ${\cal W}$ determined by the web function
$f(x,y):=(x+y)e^{-x}$, i.e.~the 3-web given by the foliations
\begin{equation}
  \label{eq:example}
  x=const, \qquad y = const, \qquad (x+y) \,e^{-x}=const.
\end{equation}
we can find in \cite{GMS} that $W$ is linearizable while
\cite{GL_cras} and \cite{GL} state the opposite.

\bigskip

In this present paper, with a computation based on the theory of
\cite{GMS}, we prove that the 3-web given by (\ref{eq:example}) is
linearizable by finding \textit{explicitly} the affine deformation
tensor.  Through this example we demonstrate the efficiency and the
correctness of the approach developed in \cite{GMS}.

\bigskip

\section{Basic notations and definitions}

Let $M$ be a $2$-dimensional differentiable manifold.
\begin{definition}
  A $3-$web on $M$ is a triple of foliations $\{F_1, F_2, F_3\}$ such
  that the tangent spaces to the leaves of any two different
  foliations are complementary subspaces of $T$.
\end{definition}
We will call the leaves of the foliations $\{F_1, F_2, F_3\}$
horizontal, vertical and transversal.  Likewise, we call their tangent
spaces {\it horizontal, vertical and transversal} and denote them by
$T^h$, $T^v$ and $T^t$. We will use Nagy's formalism (cf.
\cite{Nag}).  In particular, $h$ (resp. $v$) is the horizontal (resp.
vertical) projection, $j$ is the associated product structure, $\N$ is
the Chern connection.  By the inverse functions theorem, we can find
local coordinates ($x_1,x_2$) at a neighborhood of $p\in M$ such that
$\W$ can be written as
\begin{equation}
  \label{eq:formelocale} 
  x_1=const, \qquad x_2=const, \qquad f(x_1,x_2)=const.
\end{equation}
Using (\ref{eq:formelocale}) as a local representation of the web, at
every point the horizontal, vertical and transversal spaces are:
\begin{displaymath}
  T^h = Span\left\{ \partial_1 \right\}, \qquad T^v = Span \left\{
    \partial_2 \right\}, \qquad T^v = Span \left\{ \partial_1 - c
    (x_1,x_2) \partial_2 \right\},
\end{displaymath}
where $c(x_1,x_2):= \partial_1 f / \partial_2 f$ and $\partial_i:=
\partial/\partial x_i$.  Moreover, one has $j(\partial_1) = c
\partial_2$ and $j(\partial_2) = \frac{1}{c} \partial_1$, therefore
the basis
\begin{equation}
  \left\{e_1:= \partial_1, e_2:= c \partial_2 \right\}
\end{equation}
is an adapted basis of the web i.e. $e_1 \in T^h$, $e_2 \in T^v$ and
$j e_1 = e_2$. The Chern connection $\n$ is determined by
\begin{displaymath}
  \nabla_{\displaystyle \partial_1 } \partial_1 = \Gamma_1 \,
  \partial_1, \qquad
  \nabla_{\displaystyle \partial_2} \partial_2 = \Gamma_2 \,
  \partial_2, \qquad
  \nabla_{\displaystyle \partial_1} \partial_2 = 0, \qquad
  \nabla_{\displaystyle \partial_2} \partial_1 = 0,
\end{displaymath}
where $\displaystyle \Gamma_1= \frac{c_1}{c}$ and $\displaystyle
\Gamma_2= - \frac{c_2}{c}$.  The curvature tensor $R^\n$ of the Chern
connection is characterized by the function
\begin{displaymath}
  r:=\frac{c_1c_2-c_{12} c}{c^2} = \frac{f_{11}f_2^2f_{12} -
    f_1^2f_{12}f_{22} - f_1f_{112}f_2^2 +
    f_1^2f_{122}f_2}{f_2^2f_1^2},
\end{displaymath}
since for $i=1,2$ one has
\begin{alignat*}{1}
  R^\n \left( \partial_1, \partial_2 \right)\partial_i = r \,
  \partial_i.
\end{alignat*}

\bigskip

\section{The PDE system of linearization}

\begin{definition}
  A $3-$web on a $2$-dimensional affine space is called linear
  (resp.~parallel) if the leaves of the three foliations are straight
  lines (resp.~parallel straight lines). A $3$-web on $M$ is called
  linearizable (resp.~parallelizable) at $p\in M$ if it is equivalent
  to a linear (resp.~parallel) 3-web modulo a local diffeomorphisms.
\end{definition}
The problem of linearizability of webs can be formulated as follows:
find a torsion-free flat connection $\n^L$ such that the foliations of
the web are geodesic with respect to this connection \cite{Gol2}. The
existence of such connection is equivalent to the existence of a
symmetric (1,2)-tensor field $L$, the \textit{linearization} or
\textit{affine deformation} tensor, which satisfies to the condition
that the connection $\nabla^L$ defined as
\begin{displaymath}
  \N^L_XY := \N_XY + L(X,Y)
\end{displaymath}
preserves the web and flat.  A tensor field $L$ in $S^2T^* \otimes T$
is a linearization if and only if
\begin{alignat*}{1}
  1. & \quad vL(hX,hY) = 0,
  \\
  2. & \quad hL(vX,vY) = 0,
  \\
  3. & \quad L(hX,hY)+jL(jhX,jhY)-hL(jhX,hY)
  \\
  & \qquad - hL(hX,jhY)-jvL(jhX,hY)-jvL(hX,jhY) = 0,
  \\
  4. & \quad \N_XL (Y,Z) - \N_YL (X,Z) + L(X,L(Y,Z)) - L(Y,L(X,Z)) +
  R^\n (X,Y)Z = 0,
\end{alignat*}
holds for any $X,Y,Z \in T$, \cite{GMS}.  Using local coordinate
system, a symmetrical tensor $L = L^k_{ij}dx^i \otimes dx^j \otimes
\partial_k \in S^2T^* \otimes T$ is a linearization if and only if its
components satisfy (\ref{eq:2}) and (\ref{eq:3}), where
\begin{equation}
  \label{eq:2}
  L_{11}^2 =0, \quad L_{22}^1 =0, \quad L_{12}^2=\frac{1}{2}(L_{11}^1
  + c L_{22}^2-2cL_{12}^1),
\end{equation}
is a system of algebraic equations and 
\begin{equation}
  \label{eq:3}
  \left.
    \begin{aligned}
      r + \frac{\partial L^1_{12}}{\partial x} - \frac{\partial
        L^1_{11}}{\partial y} +L^2_{12}L^1_{12} &=0, \quad
      \\
      \frac{\partial L^2_{12}}{\partial x} - \Gamma_1L^2_{12}
      +L^2_{12}L^2_{12}- L^1_{11}L^2_{12}&=0,
      \\
      -\frac{\partial L^1_{12}}{\partial y} + \Gamma_2 L^1_{12}
      +L^2_{22} L^1_{12} - L^1_{12} L^1_{12} &=0,
      \\
      r +\frac{\partial L^2_{22}}{\partial x} -\frac{\partial
        L^2_{12}}{\partial y} - L^1_{12} L^2_{12} &=0.
    \end{aligned}
  \right\}
\end{equation}
is a system of first order quasi-linear partial differential
equations.

A tensor $L$ in $S^2T^* \otimes T$ satisfying the algebraic conditions
(\ref{eq:2}) called \textit{prelinearization}.  The prelinearizations
forms a 3 dimensional subbundle of $S^2T^* \otimes T$ which will be
denoted by ${\mathcal E}$. A section of ${\mathcal E}$ is a
linearization if it satisfies (\ref{eq:3}).

\section{Sketch of the solution of the system of linearization}

In this section we describe the steps needed to solve the PDE system
of linearization of the affine deformation tensor. It is composed of
the equations of (\ref{eq:2}) and (\ref{eq:3}).  The method used here
is the same as the one used in \cite{GMS} but here the system is
written in terms of functions and and partial derivatives instead of
tensors (see \cite{GMS}, p.~2648) and covariant derivatives.

As $\mathcal E$ is a rank-3 vector bundle, it can be parameterized by
$\{s, \, t,\,z\}$, where
\begin{displaymath}
  s := 2c L^1_{12}-cL^2_{22}, \qquad t:=L^2_{12}, \qquad z:=L^1_{12}.
\end{displaymath}
The parameter $s$ is called the {\it base of the prelinearization} and
it is a projective invariant of the linearizations: two
prelinearizations are projectively equivalent if and only if they have
the same base \cite{GMS}. Writing the system (\ref{eq:3}) with $\{s,
\, t,\,z\}$, the partial derivatives $t_1,t_2,z_1,z_1$ can be
expressed explicitly as
\begin{equation}
  \label{eq:frob}
  \begin{aligned}
    t_1 & =ts + \frac{t(f_{11}f_2-f_1f_{12} + tf_2f_1)}{f_2f_1},
    \\
    t_2 & = tz
    + \frac{f_2^3f_1 s_1 - f_2^3f_{11}s + f_{12} f_2^2f_1s
      -2f_2^2f_1^2 s_2}{3 f_2 ^2f_1^2}
    \\
    & \qquad \qquad + \frac{ f_{11}f_2^2f_{12} - f_1^2f_{12} f_{22} +
      f_1^2 f_{122} f_2 - f_1 f_{112}f_2^2 }{3 f_2 ^2f_1^2},
    \\
    z_1 & = tz+
    \frac{-s_2f_2^2f_1^2+2f_2^3f_1 s_1-2f_2^3f_{11}s+2f_{12}f_2^2f_1s
    }{3 f_2^2f_1^2}
    \\
    & \qquad \qquad +\frac{f_1f_{112}f_2^2 -f_{11}f_2^2 f_{12}
      +f_1^2f_{12}f_{22} - f_1^2 f_{122}f_2}{3 f_2^2f_1^2},
    \\
    z_2 & = z^2 - \frac{z(f_{12}f_2-f_1f_{22}+f_2^2 s)}{f_2f_1}.
  \end{aligned}
\end{equation}
By the consideration the integrability conditions $ t_{12} = t_{21}$
and $z_{12} = z_{21}$ one can realize, that the functions $t$, $z$ and
their derivatives can be eliminated. That way one obtains two second
order PDE on $s$:
\begin{alignat}{2}
  \label{eq:t}
  \mathrm{I.}) \qquad & 0= s_{11}-2c \, s_{12} && + {lower \ order \
    terms}...
  \\
  \label{eq:z}
  \mathrm{II.}) \qquad & 0 =s_{22}-\frac{2}{c}s_{12} &&+ {lower \
    order \ terms}...
\end{alignat}
There is no integrability condition coming for the first prolongation
of (\ref{eq:t}) and (\ref{eq:z}), but there is one integrability
condition coming for the second prolongation. Indeed, using the second
prolongation of (\ref{eq:t}) and (\ref{eq:z}), the equation
\begin{equation}
  \label{eq:s12}
  0 = c \, (\partial_{11} \mathrm{II} - \partial_{22} \mathrm{I})+
  2(\partial_{12} \mathrm{I} - c^2\, \partial_{12} \mathrm{II})
\end{equation}
does not contain 4th order derivatives of $s$. We express the third
order derivatives of $s$ form (\ref{eq:t}) and (\ref{eq:z}) and
substitute them into (\ref{eq:s12}), so we obtain a new equation:
\begin{equation}
  \label{eq:III}
  \mathrm{III}.) \qquad 0 = 24 \, c \, r \, s_{12} + \ {lower \ order
    \ terms}...
\end{equation}
We remark that the integrability condition (\ref{eq:III}) is
identically satisfied if $r =0$, i.e.~the web is parallelizable. If
$r  \neq 0$, the we have to push forward the computation. 

Let us suppose that $r \neq 0$. We have to consider the system formed
by (\ref{eq:z}), (\ref{eq:t}) and (\ref{eq:III}), which are second
order PDE equations on $s$.  The prolongation of these equations leads
us to 2 integrability conditions. Indeed, considering the combinations
\begin{alignat}{2}
  \label{eq:6}
  0 & = 24\, c \, r \, \mathrm{I}_2 - \mathrm{III}_1 +2 \, c \,
  \mathrm{III}_2
  \\
  \label{eq:7}
  0 & = 24\, c \, r\, \mathrm{II}_1 - \mathrm{III}_2 + \frac{2}{c}
  \mathrm{III}_1
\end{alignat}
the new equations do not contain 3rd order derivatives. Moreover the
second order derivatives of $s$ can be expressed form (\ref{eq:t}),
(\ref{eq:z}) and (\ref{eq:s12}), we can substitute them into
(\ref{eq:6}) and (\ref{eq:7}).  That way we obtain two new equations
having special forms
\begin{alignat}{2}
  \label{eq:5}  
  0 & = - 24 r \ (s_1)^2 + 48 r \ s_1 s_2 + \alpha_1 \, s_1 + \beta_1
  \ s_2 + \gamma_1
  \\
  \label{eq:8}
  0 & = c\, 24 \, r \ (s_2)^2 + 48 r \ s_1 s_2 + \alpha_2 \ s_1 +
  \beta_2 \ s_2 + \gamma_2
\end{alignat}
where $\alpha^i$, $\beta^i$ and $\gamma^i$ ($i=1,2$) are determined by
the curvature and its derivatives.  At the final step one has to take
the derivatives of (\ref{eq:5}) and (\ref{eq:8}) with respect to the
variables $x_1$ and $x_2$.  By expressing $(s_1)^2$ and $(s_2)^2$ from
(\ref{eq:5}) and (\ref{eq:8}), respectively, and by putting the
correspondent values into the prolongated derived system we obtain
a system of for equations
\begin{equation}
  \label{abc}
  \hphantom{\qquad i=1,..,4} a^is_1+b^is_2+c^i \, s_1s_2=d^i, \qquad
  i=1,..,4.
\end{equation}
where $a^i$, $b^i$, $c^i$ and $d^i$, $i=1,..,4$ are determined by the
curvature and its derivatives.  (\ref{abc}) can be considered as a
linear system in $s_1$, $s_2$ and $s_1s_2$. This system is compatible,
and the $3^{rd}$-order minors are non zero polynomials in $s$ of degree
$7$ (\cite{GMS}, p.~2652).  So there exists an open ${\mathcal U}
\subset \Bbb C^2$ on which
\begin{displaymath}
  D(s) := \left|
    \begin{array}{lll}
      a^1 & b^1 & c^1
      \\
      a^2 & b^2 & c^2
      \\
      a^3 & b^3 & c^3
    \end{array}
  \right| \neq 0.
\end{displaymath}
Solving on ${\mathcal U}$ the linear system (\ref{abc}) by the Cramer
formulas, we get:
\begin{equation}
  \label{s1s2}
  s_1 =\frac{A(s)}{D(s)}, \qquad s_2= \frac{B(s)}{D(s)},\qquad s_1s_2
  = \frac{C(s)}{D(s)},
\end{equation}
where $A$, $B$, and $C$ are given by the corresponding determinant.
These functions are polynomial in $s$. Moreover,
\begin{enumerate}
\item [(a)] using the identity $s_1 \cdot s_2 = s_1 s_2 $ and the
  expression of the corresponding terms given by (\ref{s1s2}) we
  obtain that the solution $s$ of the linearization system has to take
  his values on algebraic manifold defined by $Q_1(s)=0$, where
  \begin{math}
    Q_1(s) : = AB -CD
  \end{math}
  is polynomial in $s$ of degree 18.
\item[(b)] For the system (\ref{s1s2}) the compatibility condition is
  given by $\partial_1 s_2 -\partial_2 s_1=0$.  Using $A$, $B$ and
  $D$, we obtain that $s$ has to take its values in the algebraic
  manifold $Q_2(s)=0$ defined by this compatibility condition.
\item[(c)] Just like the first derivatives are computed in
  (\ref{s1s2}), the second order derivatives can also be expressed in
  a similar way, and using their expressions in the equations
  (\ref{eq:t}), (\ref{eq:z}), (\ref{eq:III}), (\ref{eq:5}) and
  (\ref{eq:8}) we get 5 polynomial equations in $s$: $Q_i = 0, \
  (i=3,...,7)$.
\end{enumerate}
It follows that $s=s(x_1, x_2)$ has to take its values in the
algebraic manifold ${\mathcal A} \subset E$, where
\begin{displaymath}
  \mathcal{A} : = \bigl\{Q_i=0 \ | \ i=1,..., 7 \bigl\}.
\end{displaymath}

\section{Example}

In this section we consider the 3-web determined by the web function
$f(x_1,x_2):=(x_1+x_2)e^{-x_1}$, that is the 3-web ${\cal W}$ given by
the foliations
\begin{equation}
  \label{eq:example_1}
  x_1=const, \qquad x_2 = const, \qquad (x_1+x_2)e^{-x_1}=const.
\end{equation}
The linearizability of this example was examined by Grifone, Muzsnay
and Saab in \cite{GMS} (page 2563), and the authors claimed that this
particular web is linearizable.  However, in \cite{GL} (page 38) and
\cite{GL_cras} (page 171) the authors stated the opposite.

\medskip

Let us examine this example more closely. The Chern connection of the
web $\cal W$ is determined by:
\begin{alignat*}{2}
  \nabla_{\displaystyle \partial_1 } \partial_1 & =
  \frac{1}{x_1+x_2-1} \, \partial_1, \qquad \qquad \qquad &
  \nabla_{\displaystyle \partial_1} \partial_2 &= 0,
  \\
  \nabla_{\displaystyle \partial_2} \partial_2 &= \frac{1}{1-x_1-x_2}
  \, \partial_2, \qquad \qquad \qquad & \nabla_{\displaystyle
    \partial_2} \partial_1 &= 0.
\end{alignat*}
The curvature is given by
\begin{displaymath}
  R^\n \left( \partial_1, \partial_2 \right)\partial_i =
  \frac{1}{(x_1+x_2-1)^2} \, \partial_i,
\end{displaymath}
for $i=1,2$. Therefore the Chern connection is non flat and the web
$\W$ is not parallelizable. Following the computation described in the
previous chapter one can find, that
\begin{equation}
  \label{eq:s=-1}
  s (x_1, x_2)\equiv -1
\end{equation}
is a solution for all the polynomials $Q_i(s)$, $i=1,...,7$. This
shows that the web is linearizable. Let us go further and find the
linearization explicitly. By substituting $s(x_1, x_2)\equiv -1$ into
(\ref{eq:frob}) one obtains that
\begin{equation}
  \label{eq:frob_ex}
  \begin{aligned}
    t_1 & = t^2 - t + \frac{t}{x_1+x_2-1},
    \\
    t_2 & = t z, \hphantom{\frac{1}{1}}
    \\
    z_1 & = tz - \frac{1}{(x_1+x_2-1)^2 },
    \\
    z_2 & = z^2 - \frac{2z}{x_1+x_2-1}.
  \end{aligned}
\end{equation}
There are two solutions of the differential system (\ref{eq:frob_ex}):
\begin{alignat}{1}
  \label{eq:sol_1}
  \mathrm{Solution \ 1.}& \qquad \left\{
    \begin{aligned}
      t(x_1, x_2) & = 0,
      \\
      z(x_1, x_2) & = \frac{1-x_1 -\a }{(-1+x_1+x_2)(x_2-\a )},
    \end{aligned}
  \right.
  \\
  \notag
  \\
  \label{eq:sol_2}
  \mathrm{Solution \ 2.}& \qquad \left\{
    \begin{aligned}
      t(x_1, x_2) & = \frac{(-1+x_1+x_2) \, e^{-x_1}}{(x_1+x_2)
        e^{-x_1}+\a  x_2+\bb },
      \\
      \\
      z(x_1, x_2) & = \frac{e^{-x_1} + \a -x_1 \a  + \bb }
      {\bigl((x_1+x_2) \, e^{-x_1}+\a  x_2 + \bb  \bigl) (x_1+x_2 -1)}
    \end{aligned}
  \right.
\end{alignat}
where $a$ and $b$ are arbitrary constants.

\bigskip
\bigskip

\textbf{Solution 1.} Here we consider the solution (\ref{eq:sol_1}) of
(\ref{eq:frob_ex}). Rewriting the expression of $t(x_1, x_2)$ and
$z(x_1, x_2)$ with the help of (\ref{eq:s=-1}) we can determine the
components of the affine deformation tensor $L$:
\begin{alignat*}{2}
  L_{11}^1&=-1, \qquad\qquad & L_{22}^2&= -\frac{x_2-2+2x_1+\a }
  {(x_1+x_2-1) (x_2-\a )},
  \\
  L_{12}^2&= 0, \qquad \qquad& L_{12}^1&=\frac{1-x_1
    -\a }{(-1+x_1+x_2) (x_2-\a )}.
\end{alignat*}
The deformed connection $\n^L$ in the standard base is given by the
following equations:
\begin{alignat}{1}
  \label{eq:sol_1_1}
  \n^L_{\partial_1} \partial_1 & = \n_{\partial_1} \partial_1
  +L(\partial_1, \partial_1) = \frac{c_1}{c} \partial_1+L_{11}^1
  \partial_1 = \frac{x_1+x_2-2}{1-x_1-x_2}\, \partial_1
  \\
  \label{eq:sol_1_2}
  \n^L_{\partial_1} \partial_2 & = \n_{\partial_1} \partial_2
  +L(\partial_1, \partial_2)
  = L_{12}^1 \partial_1+ L_{12}^2 \partial_2 = \frac{1-x_1
    -\a }{(-1+x_1+x_2)(x_2-\a )} \ \partial_1
  \\
  \label{eq:sol_1_3}
  \n^L_{\partial_2} \partial_1 & = \n_{\partial_2} \partial_1
  +L(\partial_2, \partial_1)
  = L_{12}^1 \partial_1+ L_{12}^2 \partial_2 = \frac{1-x_1
    -\a }{(-1+x_1+x_2)(x_2-\a )} \ \partial_1
  \\
  \label{eq:sol_1_4}
  \n^L_{\partial_2} \partial_2 & = \n_{\partial_2} \partial_2
  +L(\partial_2, \partial_2) = -\frac{c_2}{c} \partial_2+L_{22}^2
  \partial_2 = \frac{2}{\a -x_2} \ \partial_2
\end{alignat}
It is obvious that $\n^L_{\partial_i} \partial_j - \n^L_{\partial_j}
\partial_i =0$ and therefore the torsion of $\n^L$ is zero.  Moreover,
at every point the horizontal, vertical and transversal spaces are:
\begin{displaymath}
  T^h = Span\left\{ \partial_1 \right\}, \qquad T^v = Span \left\{
    \partial_2 \right\}, \qquad T^v = Span \left\{ \partial_1 - c \,
    \partial_2 \right\},
\end{displaymath}
where $c=\partial_1 f /\partial_2 f= 1-x_2-x_2$. The equation
(\ref{eq:sol_1_1}) (resp.~(\ref{eq:sol_1_4})) shows that the covariant
derivative of a horizontal (resp.~vertical) vectorfield with respect
to a horizontal (resp.~vertical) vectorfield is horizontal
(resp.~vertical). Moreover, we have
\begin{alignat*}{1}
  \n_{(\partial_1 - c\partial_2) } (\partial_1 - c\partial_2) %& = 
  & = \n_{\partial_1} \partial_1 - c\, \n_{\partial_2 } \partial_1 - c
  \n_{\partial_1 } \partial_2 - (\partial_1 c) \, \partial_2 + c
  (\partial_2 c) \, \partial_2 + c^2 \n_{\partial_2 } \partial_2
  \\
  & = \frac{2 x_1^2+x_2^2 + 3 x_1 x_2 + (a-4) (x_1+x_2) +2
  }{(-1+x_1+x_2) (a-x_2)}\, (\partial_1 - c\partial_2)
\end{alignat*}
which shows that the covariant derivative of a transversal vectorfield
with respect to a transversal vectorfield is transversal.  Direct
calculation shows that $\n^L$ is flat, that is its curvature tensor is
identically zero.

\bigskip

\textbf{Solution 2.} Here we consider the solution (\ref{eq:sol_2}) of
(\ref{eq:frob_ex}). Completing the expression of $t(x_1, x_2)$ and
$z(x_1, x_2)$ with (\ref{eq:s=-1}) we can find that the components of
$L$, the affine deformation tensor are:
\begin{alignat*}{2}
  L_{11}^1&=\frac{ (x_1 +x_2 -2)e^{-x_1} -\a x_2-\bb } {
    (x_1+x_2)e^{-x_1}+\a x_2+\bb }
  \\
  L_{22}^2&= \frac{(2-x_1-x_2)e^{-x_1} -\a (2x_1+x_2-2) +\bb }{
    (x_1+x_2-1)((x_1+x_2)e^{-x_1}+\a x_2+\bb )}
  \\
  L_{12}^1&= \frac{e^{-x_1} - \a x_1 + \a + \bb } { (x_1+x_2
    -1)\bigl((x_1+x_2) \, e^{-x_1}+\a x_2 + \bb \bigl)}
  \\
  L_{12}^2&= \frac{(x_1+x_2-1) \, e^{-x_1}}{(x_1+x_2) e^{-x_1}+\a
    x_2+\bb }
\end{alignat*}
The deformed connection $\n^L$ in the standard base is given by the
following equations:
\begin{alignat}{1}
  \label{eq:sol_2_1}
  \n^L_{\partial_1} \partial_1 & = \left( \frac{1}{x_1+x_2-1}+\frac{
      (x_1 +x_2 -2)e^{-x_1} -\a x_2-\bb } { (x_1+x_2)e^{-x_1}+\a
      x_2+\bb } \right) \, \partial_1
  \\
  \label{eq:sol_2_2}
  \n^L_{\partial_1} \partial_2 & = \frac{e^{-x_1} + \a -x_1 \a + \bb }
  {\bigl((x_1+x_2) \, e^{-x_1}+\a x_2 + \bb \bigl) (x_1+x_2 -1)} \
  \partial_1+ \frac{(x_1+x_2-1) \, e^{-x_1}}{(x_1+x_2) e^{-x_1}+\a
    x_2+\bb } \ \partial_2 ,
  \\
  \label{eq:sol_2_3}
  \n^L_{\partial_2} \partial_1 & = \frac{e^{-x_1} + \a -x_1 \a + \bb }
  {\bigl((x_1+x_2) \, e^{-x_1}+\a x_2 + \bb \bigl) (x_1+x_2 -1)} \
  \partial_1+ \frac{(x_1+x_2-1) \, e^{-x_1}}{(x_1+x_2) e^{-x_1}+\a
    x_2+\bb } \ \partial_2 ,
  \\
  \label{eq:sol_2_4}
  \n^L_{\partial_2} \partial_2 & =\frac{ -2(e^{-x_1} +\a) }{
    (x_1+x_2)e^{-x_1}+\a x_2+\bb } \ \partial_2
\end{alignat}
As in the previous case $\n^L_{\partial_i} \partial_j -
\n^L_{\partial_j} \partial_i =0$ and the torsion of $\n^L$ is zero.
Equation (\ref{eq:sol_2_1}) (resp.~(\ref{eq:sol_2_4})) shows that the
covariant derivative of a horizontal (resp.~vertical) vectorfield with
respect to a horizontal (resp.~vertical) vectorfield is horizontal
(resp.~vertical). We have
\begin{alignat*}{1}
  \n_{(\partial_1 - c\partial_2) } & (\partial_1 - c\partial_2)=
  \\
  & = \frac{(x_1+x_2)^2 e^{-x_1} +(4a+b)(x_1+x_2) -a(2x_1^2 +x_2^2
    +3x_2x_1 +2)}{ (x_1+x_2-1)( (x_1+x_2)e^{-x_1} +a x_2+b) }\,
  (\partial_1 - c\partial_2)
\end{alignat*}
which shows that the covariant derivative of a transversal vectorfield
with respect to a transversal vectorfield is transversal.  As in the
previous case, $\n^L$ is flat i.e.~its curvature tensor is identically
zero.

\bigskip \bigskip

\noindent
As the direct calculations show in both cases
\begin{enumerate}
\item the connection $\n^L$ preserves the web, that is the three
  families of leaves are auto-parallels curves with respect to it;
\item $\n^L$ is torsion free;
\item $\n^L$ is flat, that is its curvature tensor is identically
  zero.
\end{enumerate}
The properties 1.)~-~3.)~show that the corresponding affine deformation
tensor $L$ in both cases is a linearization of the web ${\cal W}$.

\bigskip

\textbf{Remark.} Solution 1.~and Solution 2.~correspond to different
linearizations. However, these linearizations are projectively
equivalent.  Indeed, the parameter $s$, called the base of the
linearization, is a projective invariant of the linearizations: two
linearizations are projectively equivalent if and only if they have
the same base.  Here the two linearizations have the same base
($s(x_1, x_2) \equiv -1$) which shows that they are projectively
equivalent.

\bigskip

\nn 
\begin{list}{}{}
\item[Zolt\'an Muzsnay] : University of Debrecen, Department of
  Mathematics, Debrecen, H-4032 PBox 12, Hungary,
\item 
  \textit{E-mail address}:  \textit{muzsnay@math.klte.hu}
\end{list}


\bigskip



\begin{thebibliography}{00}
  
\bibitem {AS} {\sc Akivis, M.A. and Shelkhov, A.M.}; {\it Geometry of
    Algebra of Multidimensional Three-Webs.} Kulwer Academic
  Publishers, Dordrecht (1992).
  
\bibitem {BB} {\sc Blaschke, W.}; {\it Einf\"uhrung in die Geometrie der
    waben.}  Birkhauser-Verlag, Basel-Stuttgart, 1955.

\bibitem{Bol1} {\sc Bol, G.}; {\it Geradlinige Kurvengewebe}. Abh.
  Math.  Sem. Univ.  Hambourg, 8, (1930), 264-270.
  
\bibitem{Bol2} {\sc Bol, G.}; {\it \"Uber Geradengewebe}. Ann. Mat.
  Pura e Appl. (4) 17 (1938), 45-58.
  
  
\bibitem{Gol2}{\sc Goldberg, V.V.}; {\it On a linearizability
    condition for a three-web on a two-dimensional manifold.}
  Differential Geometry, (Peniscola, 1988), 223-239, Lecture Notes in
  Math. 1410, Springer-Verlag, Berlin New-York, (1989).


\bibitem{GL}{\sc Goldberg,V.V.; Lychagin, V.V.}; \textit{On the
    Blaschke conjecture for 3-webs}, arXiv: math.DG/0411460

\bibitem{GL_cras}{\sc Goldberg,V.V.; Lychagin, V.V.}; \textit{On
    linearization of planar three-webs and Blaschke's conjecture},
  C.R.Acad. Sci. Paris, Ser. I.  vol. 341. num 3 (2005)

  
\bibitem{GM}{\sc Grifone, J; Muzsnay, Z.}; \textit{Variational
  Principles For Second-Order Differential Equations}, World
  Scientific, Singapore, (2000).
  
\bibitem{GMS}{\sc Grifone, J; Muzsnay, Z; Saab J.}; \textit{On the
    linearizability of 3-webs}, Nonlinear analysis 47, (2001) 2643-2654.


\bibitem{GMS_arxiv}{\sc Grifone, J; Muzsnay, Z; Saab J.};
  \textit{Linearizable 3-webs and the Gronwall conjecture}, arXiv:
  math.DG/0602535
 

\bibitem{Gro}{\sc Gronwall, T.H.} {\it Sur les \'equations entre trois
    variables repr\'esentables par des nomogrammes \`a points
    align\'es.} J. de Liouville, 8, 59-102.

    
\bibitem{Nag}{\sc Nagy, P.T.}; {\it Invariant tensorfields and the
    canonical connection of a 3-web.} Aequationes Math, 35 (1988)
  31-44.
\end{thebibliography}
\end{document}